\documentclass[12pt]{amsart}

\usepackage{amssymb}

\newtheorem{theorem}{Theorem}[section]
\newtheorem{lemma}[theorem]{Lemma}
\newtheorem{proposition}[theorem]{Proposition}
\newtheorem{remark}[theorem]{Remark}
\newtheorem{definition}[theorem]{Definition}

\numberwithin{equation}{section}

\newcommand{\too}{\longrightarrow}

\newcommand{\surj}{\twoheadrightarrow}
\newcommand{\x}{\times}
\newcommand{\ox}{\otimes}

\newcommand{\SL}{\text{SL}}

\newcommand{\GL}{\text{GL}}
\newcommand{\Pic}{\text{Pic}}

\newcommand{\tr}{\text{trace}}
\newcommand{\rk}{\text{rank}}
\newcommand{\ad}{\text{ad}}
\newcommand{\Aut}{\text{Aut}}
\newcommand{\Hom}{\text{Hom}}
\newcommand{\Id}{\text{Id}}
\newcommand{\im}{\text{image}}

\newcommand{\cD}{{\mathcal D}}
\newcommand{\cE}{{\mathcal E}}
\newcommand{\cF}{{\mathcal F}}
\newcommand{\cH}{{\mathcal H}}
\newcommand{\cM}{{\mathcal M}}

\newcommand{\cN}{{\mathcal N}}
\newcommand{\cO}{{\mathcal O}}
\newcommand{\cR}{{\mathcal R}}
\newcommand{\cS}{{\mathcal S}}
\newcommand{\cU}{{\mathcal U}}
\newcommand{\cZ}{{\mathcal Z}}

\newcommand{\CC}{{\mathbb C}}

\newcommand{\QQ}{{\mathbb Q}}
\newcommand{\RR}{{\mathbb R}}
\newcommand{\ZZ}{{\mathbb Z}}
\newcommand{\g}{\gamma}

\begin{document}

\baselineskip=15.5pt

\title[Torelli theorem for moduli of connections]{Torelli
theorem for moduli spaces of $\SL(r,\CC)$--connections on a
compact Riemann surface}

\author[I. Biswas]{Indranil Biswas}

\address{School of Mathematics, Tata Institute of Fundamental
Research, Homi Bhabha Road, Bombay 400005, India}

\email{indranil@math.tifr.res.in}

\author[V. Mu{\~n}oz]{Vicente Mu{\~n}oz}

\address{Departamento de Matem\'aticas,
Consejo Superior de Investigaciones Cient{\'\i}ficas, Serrano 113
bis, 28006 Madrid, Spain}

\curraddr{Institute for Advanced Study, Einstein Drive,
Princeton, New Jersey 08540, USA}

\email{vicente.munoz@imaff.cfmac.csic.es}

\thanks{The second author is
partially supported through MCyT grant MTM2004-07090-C03-01 (Spain)
and NSF grant No.\ DMS-0111298 (US). Any opinions,
findings and conclusions or recommendations
expressed in this material are those of the authors and do not
necessarily reflect the views of the National Science
Foundation.}



\date{}

\begin{abstract}
Let $X$ be any compact connected Riemann surface of genus $g$,
with $g\geq 3$. For any $r\,\geq\, 2$, let
$\cM_X$ denote the moduli space of
holomorphic $\SL(r,\CC)$--connections over $X$.
It is known that the biholomorphism class of
the complex variety $\cM_X$
is independent of the complex structure of $X$.
If $g\,=\,3$, then we assume that $r\,\geq\,3$.
We prove that the isomorphism class of the variety
$\cM_X$ determines the Riemann surface $X$ uniquely up to
an isomorphism.
A similar result is proved for the moduli space of
holomorphic $\GL(r,\CC)$--connections on $X$.

We also show that the Torelli theorem remains valid
for the moduli spaces of connections, as well as
those of stable
vector bundles, on geometrically irreducible smooth
projective curves defined over the field of real numbers.
\end{abstract}

\maketitle

\section{Introduction} \label{sec:intro}

Let $X$ be a compact connected Riemann surface of genus $g$,
with $g\, \geq\, 3$. Fix an integer $r\, \geq\,2$. If
$g\, =\, 3$, then we will assume that $r\, \geq\, 3$.

Let $\cM_X$ denote the moduli space parametrizing pairs of the
form $(E\, ,D)$ on $X$ satisfying the following two conditions:
\begin{itemize}
\item $E$ is a holomorphic vector bundle of rank $r$ with
determinant $\bigwedge^r E\, \cong\, \cO_X$, and

\item $D$ is a holomorphic connection on $E$ such that
the connection on $\bigwedge^r E$ induced by
$D$ coincides with the trivial connection on $\cO_X$ given
by the de Rham differential.
\end{itemize}

In other words, $\cM_X$ is the moduli space of holomorphic
${\rm SL}(r, {\mathbb C})$--connections on $X$.
This moduli space $\cM_X$ is an irreducible normal
quasiprojective variety, defined over $\CC$, of dimension
$2(r^2-1)(g-1)$.

We prove the following Torelli type theorem:

\begin{theorem}\label{thm:0}
The isomorphism class of the Riemann surface $X$ is uniquely
determined by the isomorphism class of the variety $\cM_X$. In
other words, if $Y$ is another compact connected Riemann surface
of genus $g$, and $\cM_Y$ is the moduli space of
holomorphic ${\rm SL}(r, {\mathbb C})$--connections on $Y$,
then the two varieties $\cM_X$ and $\cM_Y$ are isomorphic if
and only if the two Riemann surfaces $X$ and $Y$ are isomorphic.
\end{theorem}

Let $\widehat\cM_X$ denote the moduli space
of holomorphic $\text{GL}(r,{\mathbb
C})$--connections on $X$. This moduli space
$\widehat\cM_X$ is an irreducible normal
quasiprojective variety, defined over $\CC$, of dimension
$2(r^2(g-1)+1)$.

We prove the following analog of Theorem \ref{thm:0}:

\begin{theorem}\label{thm:00}
The isomorphism class of the Riemann surface $X$ is uniquely
determined by the isomorphism class of the variety
$\widehat\cM_X$.
\end{theorem}

The biholomorphism classes of the varieties $\cM_X$ and
$\widehat\cM_X$ depend only on $g$ and $r$. In particular,
they are independent of the complex structure of $X$
(see Proposition \ref{prop:01}).

In \cite{BM}, similar Torelli type theorems were proved
for some moduli spaces of logarithmic connections on $X$
singular over a fixed point $x_0\, \in\, X$.
More precisely, for $r\, >\, 1$ and $d$
coprime to it, it was shown that the isomorphism class of
each of the following two varieties determine
$X$ uniquely up to an isomorphism:
\begin{itemize}
\item The moduli space of logarithmic
connections $(E\, ,D)$ on $X$, where $E$ is a holomorphic
vector bundle of rank $r$ and degree $d$ with
$\bigwedge^r E \,\cong\, {\cO}_X(dx_0)$, and $D$ is a logarithmic
connection on $E$, singular exactly over $x_0$
with residue $-\frac{d}{r}\,\text{Id}_{E_{x_0}}$,
inducing the logarithmic connection on ${\cO}_X(dx_0)$
given by the de Rham differential.

\item The moduli space of logarithmic connections $(E\, ,D)$, where
$E$ is a holomorphic vector bundle of rank $r$ and degree $d$,
and $D$ is a logarithmic connection on $E$ singular exactly
over $x_0$ with residue $-\frac{d}{r}\,\text{Id}_{E_{x_0}}$.
\end{itemize}
Both these moduli spaces have the advantage of being smooth
quasiprojective varieties.
On the other hand, the moduli spaces $\cM_X$ and $\widehat{\cM}_X$
considered here are singular, which necessitates new inputs in the
proofs. The geometric significance of $\cM_X$ and $\widehat{\cM}_X$
makes them worth studying regardless of being singular.

In Section \ref{sec.real} we address Torelli type
questions for moduli spaces of objects over smooth
geometrically irreducible
projective curves defined over the field of real numbers.

\section{Preliminaries}\label{sec:1}

Let $X$ be a compact connected Riemann surface or, equivalently,
an irreducible smooth projective curve defined over the field
of complex numbers. The genus of $X$ will be denoted by $g$.

The Hodge type $(1\, ,0)$ cotangent bundle
$(T^{1,0}X)^*$ will also be denoted by ${\Omega}^{1,0}_X$, and
the holomorphic line bundle on $X$ defined by the sheaf holomorphic
sections of ${\Omega}^{1,0}_X$ will be denoted by
$K_X$. The $C^\infty$ complex line bundle $(T^{0,1}X)^*$
will be denoted by ${\Omega}^{0,1}_X$.
The trivial holomorphic line bundle
$X\x\CC$ over $X$ will also be denoted by $\cO_X$.

\begin{definition} \label{def:connection}
{\rm A \textit{holomorphic connection} on a holomorphic
vector bundle $E$ over $X$ is a first order holomorphic
differential operator
$$
\cD\, :\, E\, \too\, E\ox K_X
$$
satisfying the Leibniz identity which says that $\cD (fs) \, =\,
f\cD(s) +{\rm d}f\bigotimes s$, where $f$ (respectively, $s$) is a
locally defined holomorphic function (respectively, holomorphic
section of $E$).}
\end{definition}

Since a Riemann surface does not have nonzero $(2\, ,0)$--forms,
a holomorphic connection on a Riemann surface is automatically
flat. If $\cD$ is a holomorphic connection on $E$, then
$D+\overline{\partial}_E$ is a flat connection on $E$, where
\[
\overline{\partial}_E\, :\, E\,\too\,E\otimes {\Omega}^{0,1}_X
\]
is the Dolbeault operator for the holomorphic structure on $E$.
Conversely, if $\nabla$ is a flat connection on a $C^\infty$ vector
bundle $E$ over $X$, then the $(1\, ,0)$--component $\nabla^{1,0}$
of $\nabla$
is a holomorphic connection with respect to the holomorphic
structure on $E$ defined by the $(0\, ,1)$--component $\nabla^{0,1}$.
In particular, if a holomorphic vector bundle $V$ over $X$
admits a homomorphic connection, then $\text{degree}(V)\, =\, 0$.

\begin{remark}\label{rem.conn.}
{\rm A theorem due to Atiyah and Weil says that a holomorphic
vector bundle $E$ over $X$
admits a holomorphic connection if and only if each indecomposable
holomorphic direct summand of $E$ is of degree zero
\cite[p.\ 203, Theorem 10]{At}, \cite{We1}. Let $E$ be a holomorphic
vector bundle over $X$ that admits a holomorphic connection. If
$D$ is a holomorphic connection $E$ and
\[
\theta\, \in\, H^0(X,\, {\rm End}(E)\otimes K_X)\, ,
\]
then the differential operator $D+\theta\, :\,
E\, \too\, E\bigotimes K_X$ is also a holomorphic connection
on $E$; here $\text{End}(E)\,=\, E\bigotimes_{{\mathcal O}_X}E^*$.
This operation defines an action of
$H^0(X,\, \text{End}(E)\bigotimes K_X)$
on the space of all holomorphic connections on $E$. It
is easy to see that this way
the space of all holomorphic connections
on $E$ is an affine space for the vector
space $H^0(X,\, \text{End}(E)\bigotimes K_X)$.}

{\rm As before, let $E$ be a holomorphic
vector bundle over $X$ that admits a
holomorphic connection, and let $r$ be its rank. Fix
a holomorphic connection $D_0$ on the determinant line bundle
$\bigwedge^{r}E$. Consider the space of all holomorphic
connections on $E$ that induce the fixed connection $D_0$
on $\bigwedge^{r}E$. We will show that
this space is nonempty. For this purpose,
fix any holomorphic connection $D_E$ on $E$. Let
$D'$ be the holomorphic connection on $\bigwedge^{r}E$
induced by $D_E$. So
$$
\omega \,:=\, D'-D_0\, \in\, H^0(X,\, K_X)\, .
$$
Consider the holomorphic connection
$$
D'_E \, :=\, D_E - {\rm Id}_E\otimes \frac{\omega}{r}
$$
on $E$. The holomorphic connection on $\bigwedge^{r}E$
induced by $D'_E$ clearly coincides with $D_0$.}

{\rm Given a holomorphic connection $D$ on $E$
that induces $E_0$ on $\bigwedge^{r}E$, for any
\[
\theta\, \in\, H^0(X,\, {\rm ad}(E)\otimes K_X)\, ,
\]
where $\text{ad}(E)\, \subset\, \text{End}(E)$ is the subbundle
of corank one given by the sheaf of trace zero endomorphisms,
the holomorphic connection on $E$ given by the differential operator
$D+\theta$ has the property that the induced connection on
$\bigwedge^{r}E$ coincides with $D_0$. Conversely,
if $D$ and $D'$ are two holomorphic connections on $E$ inducing
the connection $D_0$ on $\bigwedge^{r}E$, then
\[
D'-D\, \in\, H^0(X,\, {\rm ad}(E)\otimes K_X)\, .
\]
Therefore, the space of all holomorphic connections
on $E$ that induce the fixed connection $D_0$
on $\bigwedge^{r}E$ is an affine space for the
vector space $H^0(X,\, \text{ad}(E)\bigotimes K_X)$}.$\hfill{\Box}$
\end{remark}

Given a holomorphic connection $(E\, ,D)$ on $X$, a
holomorphic subbundle $F$ of $E$ is said to be \textit{left invariant}
by $D$ if the differential operator $D$ sends any locally defined
holomorphic section of $F$ to a section of $F\bigotimes K_X$.
A holomorphic connection $(E\, ,D)$ is called \textit{reducible}
if there is a holomorphic subbundle $F$ of $E$ with $1\leq \,
\text{rank}(F)\, <\, \text{rank}(E)$ which is left
invariant by $D$.
A holomorphic connection $(E\, ,D)$ is called \textit{irreducible}
if it is not reducible.

\begin{remark}\label{rem.conn2}
{\rm Let $E$ be a holomorphic vector bundle over $X$ equipped
with a holomorphic connection $D$. If a subbundle $F\, \subset\,
E$ is left invariant by $D$, then $F$ admits a holomorphic
connection. In that case we have
\[
\text{degree}(E) \, =\, \text{degree}(F)\, =\, 0\, .
\]
Therefore, any holomorphic connection on $X$ is semistable
\cite[p.\ 88]{Si3}. Furthermore, a holomorphic connection $D$
is stable if and only if it is irreducible \cite[p.\ 88]{Si3}.
Two holomorphic connections on $X$ are called
\textit{Jordan--H\"older equivalent} if their
semisimplifications are isomorphic \cite[p.\ 90]{Si3}.
Therefore, any holomorphic connection is Jordan--H\"older equivalent
to a unique, up to an isomorphism, direct sum of irreducible
connections. A direct sum of irreducible connections is also
called a polystable connection.}$\hfill{\Box}$
\end{remark}

Let $D_0$ denote the connection on the
trivial line bundle $\cO_X$ defined by the de Rham differential
that sends any locally defined holomorphic function $f$ on
$X$ to the holomorphic
one--form ${\rm d} f$. For any integer $r \, \geq\, 1$,
let $\cM_X$ denote the moduli
space of pairs $(E\, , D)$ of the following type:
\begin{itemize}
 \item $E$ is a holomorphic vector bundle of rank $r$ over
$X$ with $\bigwedge^r E \, \cong\, \cO_X$,
 \item $D$ is a holomorphic connection on $E$, and
 \item the flat connection on $\bigwedge^r E$
 induced by $D$ has trivial monodromy.
 \end{itemize}
The last condition on $D$ is equivalent to the
following condition: the connection on $\cO_X$ given by $D$
coincides with $D_0$. We note that
using an isomorphism $\bigwedge^r E\, \longrightarrow\, \cO_X$,
a connection on $\bigwedge^r E$ gives a connection on
$\cO_X$ which is independent of the choice of the isomorphism
between $\bigwedge^r E$ and $\cO_X$.

Therefore, $\cM_X$ is the moduli space of holomorphic
$\text{SL}(r,{\mathbb C})$--connections on $X$.
See \cite{Si2}, \cite{Si3} for the construction of the moduli
space $\cM_X$. The scheme $\cM_X$ is a reduced and irreducible
normal quasiprojective variety defined over $\CC$, and its
(complex) dimension is $2(r^2-1)(g-1)$, where $g$ is the genus
of $X$ \cite[p.\ 70, Theorem 11.1]{Si3}. The closed points of
the moduli space $\cM_X$ are in bijection with all
Jordan--H\"older
equivalence classes of holomorphic connections that satisfy
the above three conditions (see Remark \ref{rem.conn2} for
the definition of the Jordan--H\"older equivalence).

Let $\widehat\cM_X$ denote the moduli space of $\text{GL}(r,
{\mathbb C})$--connections on $X$. Therefore,
$\widehat\cM_X$ parametrizes Jordan--H\"older
equivalence classes of holomorphic connections of rank $r$.
The moduli space $\widehat\cM_X$ is an irreducible normal
quasiprojective variety defined over
$\CC$, and its (complex) dimension is $2r^2(g-1)+2$
(see \cite[p.\ 70, Theorem 11.1]{Si3}).
Both $\cM_X$ and $\widehat\cM_X$ are singular varieties.

\begin{proposition} \label{prop:01}
The biholomorphism class of $\cM_X$ is independent of the complex
structure of $X$; it depends only on the integers $g$ and $r$.
The same statement holds for $\widehat\cM_X$.
\end{proposition}

\begin{proof}
There is a canonical biholomorphism
 \begin{equation}\label{eqn:extra2}
\cM_X\,\stackrel{\sim}{\too}\, \cR_{g,r}\, :=
\, \Hom (\pi_1(X)\, , \SL(r, \CC))// \SL(r, \CC) \, ,
  \end{equation}
which sends any flat connection to its monodromy
\cite[p.\ 26, Theorem 7.1]{Si3}; in $\cR_{g,r}$, two
equivalence classes of homomorphisms from
$\pi_1(X)$ to $\SL(r, \CC)$ are identified if and only if their
semisimplifications are isomorphic. The biholomorphism class of
the complex variety
$\Hom (\pi_1(X)\, , \SL(r, \CC))// \SL(r, \CC)$ depends only on
the isomorphism class of the group $\pi_1(X)$ and the integer $r$.

Similarly, the monodromy map gives a biholomorphism
$$
\widehat\cM_X\,\stackrel{\sim}{\longrightarrow}\,\widehat\cR_{g,r}
\,:=\,\Hom (\pi_1(X)\, ,\text{GL}(r, \CC))// \text{GL}(r, \CC)
$$
\cite[p.\ 26, Theorem 7.1]{Si3}; as before,
in $\widehat\cR_{g,r}$, two equivalence classes of
homomorphisms from $\pi_1(X)$ to $\text{GL}(r, \CC)$ are identified
if their semisimplifications are isomorphic. Again, the
biholomorphism class of
$\Hom (\pi_1(X)\, , \GL(r, \CC))// \GL(r, \CC)$ depends only on
the isomorphism class of the group $\pi_1(X)$ and the integer $r$.
\end{proof}

\section{The second intermediate Jacobian of the moduli
space}\label{sec3}

We continue with the notation of the previous section.
Henceforth, we will assume that $g\, =\, \text{genus}(X)\, \geq\,
3$, and $r\, \, \geq\, 2$. If $g\, =\,3$, then we will assume
that $r\, \, \geq \, 3$.

Let
\begin{equation}\label{eqn:cU}
\cU\, \subset\, \cM_X
\end{equation}
be the Zariski open subset parametrizing all
holomorphic connections $(E\, ,D)$ such that
the underlying vector bundle $E$ is stable. The openness of this
subset follows from \cite[p.\ 635, Theorem 2.8]{Ma},
\cite[p.\ 182, Proposition 10]{Sh}. Let $\cN_X$ denote the moduli
space parametrizing all isomorphism classes of stable vector
bundles $E$ over $X$ with $\rk(E)\,=\, r$ and $\bigwedge^r
E\,\cong\,\cO_X$. This moduli space $\cN_X$ is an irreducible
smooth quasiprojective variety, of dimension $(r^2-1)(g-1)$,
defined over $\CC$. Let
\begin{equation}\label{eqn:Ph}
\Phi\, :\, \cU\, \too\, \cN_X
\end{equation}
be the forgetful morphism that sends any pair $(E\, ,D)$ to $E$.

Any holomorphic vector bundle
$E\, \in\, \cN_X$ admits a unitary flat connection
\cite{NS}, hence the projection $\Phi$ in Eqn.\
\eqref{eqn:Ph} is surjective. Any holomorphic connection
on a stable vector bundle is irreducible. Hence from Remark
\ref{rem.conn.}
it follows that $\Phi$ makes $\cU$ an affine bundle over $\cN_X$.
More precisely, $\cU$ is a torsor over
$\cN_X$ for the holomorphic cotangent bundle $T^*\cN_X$. This
means that the fibers of the vector bundle $T^*\cN_X$ act freely
transitively on the fibers of $\Phi$.

As $\cM_X$ is irreducible, and $\cU$ is nonempty, the
open subset $\cU\,\subset\, \cM_X$ is Zariski dense.

\begin{lemma}\label{lem:1}
Let $\cZ\, :=\, \cM_X\setminus \cU$ be the complement of the
Zariski open dense subset $\cU$ of $\cM_X$
in Eqn.\ \eqref{eqn:cU}. The codimension of
the Zariski closed subset $\cZ$ of $\cM_X$ is at least three.
\end{lemma}

\begin{proof}
Let
\begin{equation}\label{eqn:d.Y}
Y\, \subset\, \cM_X
\end{equation}
be the Zariski closed subset of the moduli
space that parametrizes holomorphic connections that
are not stable. Therefore, $Y$ is the complement of the
stable locus in $\cM_X$ (see Remark \ref{rem.conn2}).
Given any holomorphic connection on $X$,
there is a canonically associated polystable holomorphic connection
which is unique up to an isomorphism; see Remark \ref{rem.conn2}.
Two holomorphic connections
on $X$ give the same point in the moduli space $\cM_X$ if and only
if their Jordan--H\"older equivalence classes coincide.

We will first show that the codimension of the subset $Y$ in
Eqn.\ \eqref{eqn:d.Y} is at least three. For that purpose, given
any integer $1\,\leq\, \ell\, \leq\, r$, let ${\cM}^\ell_X$
denote the moduli space of all $\text{SL}(\ell,
{\mathbb C})$--connections on $X$. So $\cM_X\, =\, {\cM}^r_X$.
The moduli space of all $\text{GL}(1,{\mathbb C})$--connections on
$X$ will be denoted by ${\widehat\cM}^1_X$.

For each $1\,\leq\, \ell\, <\, r$, we have a canonical morphism
\[
f_\ell\, :\, {\cM}^\ell_X\times {\cM}^{r-\ell}_X\times
{\widehat\cM}^{1}_X \, \longrightarrow\, \cM_X
\]
defined by
\[
((E_1\, , D_1)\, , (E_2\, , D_2)\, , (L\, , D_3))\, \longmapsto\,
((E_1\otimes L^{\otimes (r-\ell)})\oplus (E_2\otimes (L^*)^{\otimes
\ell})\, , D)\, ,
\]
where $D$ is the connection on $((E_1\bigotimes
L^{\otimes (r-\ell)})\bigoplus (E_2\bigotimes (L^*)^{\otimes \ell})$
induced by $D_1$, $D_2$ and $D_3$. Clearly,
\[
\text{image}(f_\ell) \, \subset\, Y\, ,
\]
where $Y$ is the subset in Eqn.\ \eqref{eqn:d.Y}. Furthermore,
\begin{equation}\label{eqY}
Y\, =\, \bigcup_{\ell=1}^{r-1}\text{image}(f_\ell)\, .
\end{equation}

We have $\dim {\cM}^\ell_X\, =\, 2(\ell^2-1)(g-1)$,
and $\dim {\widehat\cM}^{1}_X\, =\, 2g$,
where $g\, =\, \text{genus}(X)$. Hence
$$
\begin{aligned}
\dim \text{image}(f_\ell) & \, \leq\, 2(g-1)(r^2-2r\ell
+2\ell^2-1)+2\\
& \,\leq\, 2(g-1)(r^2-2r +1)+2 \\
& \,=\, \dim {\cM}_X  +2 -4r (g-1) \\
& \,\leq \, \dim {\cM}_X -14\, .
\end{aligned}
$$
Therefore, using Eqn.\ \eqref{eqY} we conclude that the codimension
of the Zariski closed subset $Y$ of ${\cM}_X$ is at least three.

Consider the complement $\cU'=\cM_X\setminus Y$. Then
$\cU\, \subset\, \cU'$, and $\cU$ is Zariski open in
$\cU'$. Let
 \begin{equation}\label{eqn:cZ}
 {\cZ}'\, :=\, \cU'\setminus \cU \, \subset\, \cU'
 \end{equation}
be the Zariski closed subset.
Since the codimension of $Y\, \subset\,{\cM}_X$
is at least three, to complete the proof of the lemma it suffices to
show that the codimension of ${\cZ}'$ in $\cU'$ is at least three.

We will first
show that any $(E\, ,D)\, \in\, \cU'$ is simple, i.e.,
$\dim \Aut(E\, , D)\,=\,1$.
For that purpose, consider the automorphism of the underlying vector
bundle $E$ given by an automorphism $T$ of $(E\, ,D)$, which will also
be denoted by $T$. Let $\lambda \,\in\, \CC$ be an eigenvalue
of $T(x)$ for some $x\,\in\, X$. Then $T'\, :=\, T- \lambda\cdot \Id_E$
is also an endomorphism of $(E\, , D)$. Since the endomorphism $T'$ of
$E$ fails to be an isomorphism over the above point $x$, it follows
that $F\,:=\,\text{kernel} (T')\,\neq\, 0$. Indeed, the coherent
subsheaf of $E$ generated by the parallel translations of the
subspace $\text{kernel} (T'(x))\, \subset\, E_x$ is
contained in $\text{kernel} (T')$.
Clearly $F$ is left invariant by $D$. On the other hand
$D$ is irreducible. Therefore, we have $F\, =\, E$. This
implies that $T'\,=\,0$, and hence $T\,=\,\lambda\cdot\Id_E$.

Fix a holomorphic vector bundle $E$ which admits an
irreducible holomorphic connection.
Consider the space $\cD_E$ consisting of all irreducible
holomorphic connections
$D$ on $E$ such that $(E\, , D)\, \in\, \cU'$. From the openness
of the stability condition it follows that $\cD_E$ is a Zariski
open dense subset of an affine space for the vector
space $H^0(X,\, \text{ad}(E)\bigotimes K_X)$ (see
Remark \ref{rem.conn.}). Let

\begin{equation}\label{varphi}
\varphi\, :\, \cD_E\, \longrightarrow\, \cU'
\end{equation}
be the obvious tautological map.

We have shown above that for the natural action of
the global automorphism group $\Aut(E)$ on $\cD_E$, the isotropy
at any point of $\cD_E$
is the subgroup defined by all nonzero scalar
multiplications. Therefore, for the map $\varphi$ in Eqn.\
\eqref{varphi},
\[
 \begin{aligned}
 &\dim \varphi(\cD_E) \, =\,\dim H^0(X,\, \ad(E)\ox K_X) -
 (\dim \Aut(E) - 1) \\
 &=\,
 \dim H^1(X,\, \ad(E)) - \dim H^0(X,\, \ad(E)) \, =\, (r^2-1)(g-1)\, ,
 \end{aligned}
\]
where the last equality is the Riemann--Roch formula, and the
middle equality follows from the Serre duality.
Hence, to prove that ${\cZ}'$ defined in Eqn.\ \eqref{eqn:cZ} has
codimension at least three in $\cU'$ it suffices to show
that the family of non--stable holomorphic vector bundles
admitting an irreducible holomorphic
connection has dimension no more than $(r^2-1)(g-1) -3$.

Let $\cF$ denote the family of non--stable vector bundles
admitting an irreducible connection.

First take any $E\, \in\, \cF$ which is not semistable. Let
\[
0\, =\, E_0 \, \subset\, E_1 \, \subset\, E_2
 \, \subset\, \cdots\, \subset\, E_{\ell-1}
 \, \subset\, E_{\ell} \, =\, E
\]
be the Harder--Narasimhan filtration of $E$.
We recall that the collection of pairs of integers
$\{(\text{rank}(E_i)\, ,
\text{degree}(E_i))\}_{i=1}^\ell$ is called the
\textit{Harder--Narasimhan polygon} of $E$ (see
\cite[p.\ 173]{Sh}).
The space of all isomorphism classes of holomorphic vector
bundles over $X$, whose Harder--Narasimhan polygon coincides
with that of the given vector bundle $E$, is of dimension at
most $r^2(g-1)-(r-1)(g-2)$ (this follows from
\cite[p.\ 247--248]{Bh}; see also \cite{BM}).
Therefore the locus $\cF_1$ of non--semistable vector
bundles in ${\cF}$ has dimension at most $r^2(g-1)-
(r-1)(g-2)-g$. Also, only finitely many Harder--Narasimhan
polygons occur in a bounded family of vector bundles
over $X$ \cite[p.\ 183, Proposition 11]{Sh}.
Therefore,
\[
\dim \cF_1 \,\leq\, r^2(g-1)-(r-1)(g-2)-g \, \leq\,
(r^2-1)(g-1) -3\, .
\]

Finally, take any $E\, \in\, \cF$ which is semistable but
not stable. Consider a Jordan--H\"older filtration of $E$ given by
\begin{equation}\label{eqn:j-h-filtration}
0\, =\, E_0 \, \subset\, E_1 \, \subset\, E_2
\, \subset\, \cdots\, \subset\, E_{\ell-1}
\, \subset\, E_{\ell} \, =\, E\, ,
\end{equation}
where $Q_i\, :=\, E_i/E_{i-1}$, $1\,\leq\, i \,\leq\, \ell$,
is a stable vector bundle of degree zero.
The vector bundle
\begin{equation}\label{eqn:gradation}
Q\, :=\,\bigoplus_{i=1}^\ell Q_i
 \end{equation}
is the {\em graduation} of $E$. The isomorphism class of
$Q$ depends only on $E$. The rank of $Q_i$ will be denoted
by $r_i$.

The dimension of the moduli space of all stable vector bundles
over $X$ of rank $r_i$ is $r_i^2(g-1) +1$. Therefore, the
dimension of the space of all graduations which are direct sum
of stable vector bundles of ranks $r_1\, ,r_2\, ,\cdots\, ,
r_\ell$ is
\begin{equation}\label{eq1}
\ell-g+ \sum_{i=1}^{\ell} r_i^2(g-1)
\end{equation}
(note that $\bigotimes_{i=1}^\ell \bigwedge^{r_i} Q_i\, \cong\,
{\mathcal O}_X$).

We will now calculate the dimension of the space of all
semistable vector with a fixed graduation.

Let $\cS$ be the space of all isomorphism
classes of semistable vector bundles of degree zero
whose Jordan--H\"older filtration
has the fixed graduation $Q$ in Eqn.\ \eqref{eqn:gradation}.
We may assume that isomorphism classes of
all the $Q_i$ are distinct. Indeed, if
some $Q_{i_1}$ is assumed to be isomorphic to
some $Q_{i_2}$, then the dimension of all possible
graduations of that type becomes smaller than the
number in Eqn.\ \eqref{eq1}. Using this it is easy to see
that it suffices to consider the case where all $Q_i$ are
distinct.

For $2\,\leq\, j\,\leq\, \ell$,
the dimension of the spaces of equivalence classes of extensions
$$
0\,\too\, E_{j-1}\,\too\, F\,\too\, Q_j \,\too\, 0\, ,
$$
is
$$
h^1(Q_j^*\ox E_{j-1})-1\,=\,r_j\left(\sum_{k<j}r_k\right)(g-1)-1\, ;
$$
this follows from the Riemann-Roch, the condition that
$\text{degree}(Q_i)\,=\, 0$ for all $i$, and the assumption
that all the $Q_i$ are distinct. Therefore,
\begin{equation}\label{eq2}
\dim \cS\, \leq\, -\ell+
\sum_{i=1}^\ell \sum_{k<i}r_ir_k(g-1)\, .
\end{equation}

Combining Eqn.\ \eqref{eq1} with Eqn.\ \eqref{eq2} we conclude that
the dimension of the space of all semistable vector
bundles having graduation of numerical
type $\{r_1\,, r_2\, ,\cdots\, ,r_\ell\}$
is maximum when $\ell\,=\,2$ and $\{r_1\, ,r_2\}\,=\, \{1\, ,r-1\}$.
Combining Eqn.\ \eqref{eq1} and Eqn.\ \eqref{eq2},
the dimension of all semistable vector
bundles with graduation type $\{1\, ,r-1\}$ is at most
\[
2-g +((r-1)^2+1)(g-1)+(r-1)(g-1)-1
\, =\, (r^2-r)(g-1)\, \leq\, (r^2-1)(g-1)-3\, .
\]
This completes the proof of the lemma.
\end{proof}

For any complex algebraic variety $Y$, the torsionfree part
$H^i(Y,\, \ZZ)/\text{Torsion}$ is equipped with a mixed Hodge
structure for all $i\, \geq\, 0$ \cite{De2}, \cite{De3}.

\begin{proposition}\label{prop2}
Let $\cM_X^o$ be the smooth locus of $\cM_X$. Then
the mixed Hodge structure $H^3(\cM_X^o,\,\ZZ)/{\rm Torsion}$
is isomorphic to the mixed Hodge
structure $H^3(\cN_X,\, \ZZ)/{\rm Torsion}$,
where $\cN_X$ is the moduli space in Eqn.\ \eqref{eqn:Ph}.
\end{proposition}

\begin{proof}
Consider the diagram of morphisms
 \begin{equation}\label{eqn:diag.}
 \cN_X \, \stackrel{\Phi}{\longleftarrow}\, \cU\,
 \stackrel{\iota}{\too}\, \cM_X^o\, ,
 \end{equation}
where $\Phi$ is the projection in Eqn.\ \eqref{eqn:Ph}, and
$\iota$ is the inclusion map.
Since $\cU\,\stackrel{\Phi}{\too}\, \cN_X$ is a fiber bundle with
contractible fibers, the induced homomorphism
 \begin{equation}\label{eqn:is.}
 \Phi^* \, :\, H^i(\cN_X,\, \ZZ)\, \too\,
H^i(\cU,\, \ZZ)
\end{equation}
is an isomorphism for all $i\, \geq\, 0$. Therefore, $\Phi$
induces an isomorphism of the two mixed Hodge structures
$H^i(\cN_X,\, \ZZ)/\text{Torsion}$ and
$H^i(\cU,\, \ZZ)/\text{Torsion}$.

Let
\begin{equation}\label{eqn:ex.m.}
H^3(\cM_X^o\, , \cU\, ,\ZZ)\, \too\, H^3(\cM_X^o\, ,\ZZ)\,
\stackrel{\iota^*}{\too}\, H^3(\cU\, ,\ZZ)
\,\too\, H^4(\cM_X^o\, , \cU\, ,\ZZ)
 \end{equation}
be the long exact sequence of relative cohomologies.
We note that Eqn.\ \eqref{eqn:ex.m.} is an exact
sequence of mixed Hodge structures \cite[p.\ 43, Proposition
(8.3.9)]{De3}. From Lemma \ref{lem:1} (and because
$\cM_X^o$ is smooth) we know that
 \begin{equation}\label{co.va.}
H^i(\cM_X^o\, ,\cU\, ,\ZZ)\, =\, 0
 \end{equation}
for all $i\, \leq\, 4$. Therefore, the homomorphism
$\iota^*$ in Eqn.\ \eqref{eqn:ex.m.} is an isomorphism.
Consequently, the composition homomorphism
\[
(\iota^*)^{-1}\circ \Phi^*\, :\, H^3(\cN_X,\, \ZZ)
\, \too\, H^3(\cM_X^o\, ,\ZZ)\, ,
\]
where $\Phi^*$ is constructed in Eqn.\ \eqref{eqn:is.},
is the required isomorphism in the statement of the proposition.
\end{proof}

Intermediate Jacobians for mixed Hodge structures was
introduced in \cite{Ca} (see \cite[p.\ 110]{Ca}). Let
$$
 J^2(\cM_X^o)\, :=\, H^3(\cM_X^o,\, \CC)/ (F^2H^3(\cM_X^o,\, \CC)+
 H^3(\cM_X^o,\, \ZZ))
$$
be the intermediate Jacobian of the mixed Hodge structure
$H^3(\cM_X^o)$. The intermediate Jacobian of
any mixed Hodge structure is a generalized torus \cite[p.\
111]{Ca}.

\begin{proposition}\label{prop:3}
The intermediate Jacobian $J^2(\cM_X^o)$ is isomorphic to
$J^2(\cN_X)$, which is isomorphic to the Jacobian
${\rm Pic}^0(X)$ of the Riemann surface $X$.

Given a smooth family $X_T\, \longrightarrow\,
T$ of irreducible
smooth projective curves of genus $g$, let $\cM_T$
be the corresponding relative family of moduli spaces
of holomorphic ${\rm SL}(r, {\mathbb C})$--connections, and let
$\cM_T^o$ be the family consisting of the smooth locus of these moduli
spaces.
Then the relative family of Jacobians over $T$
for the family $X_T$ is isomorphic to
the relative family of second intermediate
Jacobians for $\cM_T^o$.

\end{proposition}

\begin{proof}
There is a natural isomorphism of $J^2(\cN_X)$
with $\Pic^0(X)$ \cite[p.\ 2, Theorem 1.0.2(b)]{AS}.
Therefore, the first part of the proposition
follows from Proposition \ref{prop2}.

To prove the second part we note that the construction
of the isomorphism of $J^2(\cN_X)$ with $\Pic^0(X)$
in \cite[p.\ 2, Theorem 1.0.2(b)]{AS} works over a family
of curves. Also, the construction of the isomorphism in
Proposition \ref{prop2} evidently works over a family
of curves.
This completes the proof of the proposition.
\end{proof}

In the next section we will investigate the cohomology
ring of $\cM_X^o$.

\section{Cohomology of moduli spaces}\label{sec:polarization}

We start with a proposition.

\begin{proposition}\label{prop:H2}
 $H^2(\cM_X^o,\, \ZZ)\, =\, \ZZ$.
\end{proposition}

\begin{proof}
Consider the long exact sequence of relative cohomologies
 $$
 H^2(\cM_X^o\, , \cU\, ,\ZZ)\, \too\, H^2(\cM_X^o\, ,\ZZ)\,
\stackrel{\iota^*}{\too}\,
 H^2(\cU\, ,\ZZ)\, \too\, H^3(\cM_X^o\, , \cU\, ,\ZZ)
 $$
given by $\iota$ in Eqn.\ \eqref{eqn:diag.}. Using
Eqn.\ \eqref{co.va.}, from it we conclude that
\[
 H^2(\cM_X^o,\, \ZZ)\, =\, H^2(\cU,\, \ZZ)\, .
\]
Hence, in view of Eqn.\ \eqref{eqn:is.},
to prove the proposition it suffices to show that
\begin{equation}\label{eq.Br.}
H^2(\cN_X,\, \ZZ)\, =\, \ZZ\, .
\end{equation}

Let $\text{NS}(\cN_X)$ be the N\'eron--Severi group
of $\text{NS}(\cN_X)$. Let
\[
\text{Br}'(\cN_X)\, :=\, H^2_{\text{\'{e}t}}(\cN_X,\,
{\mathbb G}_m)\, =\, H^2_{\text{\'{e}t}}
(\cN_X,\, {\mathcal O}_{\cN_X}^*)
\]
be the cohomological Brauer group. It is known that
\[
\text{NS}(\cN_X)\, =\, \text{Pic}(\cN_X)\, =\,
\ZZ
\]
(see \cite[p.\ 55, Th\'eor\`eme B]{DN}). Therefore, using
the exact sequence in \cite[p.\ 145, (8.7)]{Gr}
it can be deduced that Eqn.\ \eqref{eq.Br.}
holds provided $\text{Br}'(\cN_X)$ is a finite group. More
precisely, from \cite[p.\ 145, (8.7)]{Gr} we know that
if $\text{Br}'(\cN_X)$ is a finite group, then
$$
H^2(\cN_X,\, {\mathbb Z}/p{\mathbb Z})\, =\,
\text{NS}(\cN_X)\otimes_{\mathbb Z} {\mathbb Z}/p{\mathbb Z}
\, =\, {\mathbb Z}/p{\mathbb Z}
$$
for all but finitely many primes $p$. If $H^2(\cN_X,\, {\mathbb
Z}/p{\mathbb Z})\, =\, \text{NS}(\cN_X)\bigotimes_{\mathbb Z}
{\mathbb Z}/p{\mathbb Z} \, =\, {\mathbb Z}/p{\mathbb Z}$ for all but
finitely many primes $p$, then it is straight--forward to deduce
that
\begin{equation} \label{eqn:H2-Pic}
 H^2(\cN_X,\, {\mathbb Z})\, =\, \text{NS}(\cN_X) \, =\, {\mathbb
 Z}\, .
\end{equation}
Therefore, Eqn.\ \eqref{eq.Br.} holds provided
$\text{Br}'(\cN_X)$ is a finite group.

We will complete the proof of the proposition by showing that
$\text{Br}'(\cN_X)$ is a finite group.

Although there is no universal vector bundle over
$X\times \cN_X$, there is a canonical
universal projective bundle
\[
f'\, :\, \widehat{\mathbb P}\, \longrightarrow\, X\times\cN_X
\]
(see \cite[p.\ 6, Theorem 2.7]{BBNN}).
Fix a point $x_0\, \in\, X$. Let
\begin{equation}\label{f}
f\, :\, \widehat{\mathbb P}_{x_0}\, \longrightarrow\,
\{x_0\}\times \cN_X\, =\, \cN_X
\end{equation}
be the restriction of $f'$.
There is a Zariski open dense subset
\begin{equation}\label{f1}
{\mathcal U}'\, \subset\, \widehat{\mathbb P}_{x_0}
\end{equation}
such that there is a universal vector bundle
\begin{equation}\label{evb}
{\mathcal E}\, \longrightarrow\, X\times {\mathcal U}'
\, \subset\, X\times \widehat{\mathbb P}_{x_0}
\end{equation}
satisfying the condition that
for each point $t\, \in\, {\mathcal U}'$,
the restriction of ${\mathcal E}$ to $X\times \{t\}$ is
in the isomorphism class of stable vector bundles represented
by the point $f(t)\,\in\, \cN_X$. To construct
${\mathcal E}$, consider the rational map
\[
\beta\, :\, \widehat{\mathbb P}_{x_0}\, \dasharrow\,
{\mathcal N}_X(-x_0)
\]
defined by Hecke transformation, where ${\mathcal N}_X(-x_0)$ is
the moduli space of stable vector bundles $V$ over $X$ of rank $r$
with $\bigwedge^r V\, =\, {\mathcal O}_X(-x_0)$. The open subset
${\mathcal U}'$ in Eqn.\ \eqref{f1} is contained in the domain of
$\beta$. Let ${\mathcal V}$ denote the universal vector bundle over
$X\times {\mathcal N}_X(-x_0)$. The vector bundle ${\mathcal E}$
in Eqn.\ \eqref{evb} is obtained using the tautological Hecke
transformation on $(\text{Id}_X\times \beta)^*{\mathcal V}$.

Let
\begin{equation}\label{f.st}
f^*\, :\, \text{Br}'(\cN_X)\,\longrightarrow\,
\text{Br}'(\widehat{\mathbb P}_{x_0})
\end{equation}
be the pull back homomorphism for the projection $f$ in
Eqn.\ \eqref{f}. The relative Picard group
$\text{Pic}(\widehat{\mathbb P}_{x_0}/\cN_X)$ for the
projection $f$ in Eqn.\ \eqref{f}
is isomorphic to $\mathbb Z$.
The relative anticanonical line bundle for the projection $f$
is nonzero in $\text{Pic}(\widehat{\mathbb P}_{x_0}/\cN_X)$,
and it restricts to ${\mathcal O}(r)$ on each fiber of $f$.
Therefore, from the exact sequence in \cite[p.\ 127, (5.9)]{Gr}
we conclude that the kernel of the homomorphism $f^*$
in Eqn.\ \eqref{f.st} is a finite group. Consequently, in
order to prove that $\text{Br}'(\cN_X)$ is a finite
group it suffices to show that
\begin{equation}\label{f.st2}
f^*\, =\, 0\, .
\end{equation}

There is a Zariski open subset
\begin{equation}\label{eqs}
{\mathbb P}^s_0\, \subset\, \overline{{\mathbb P}^s_0}
\end{equation}
of a complex projective space (following the notation of \cite{DN},
the superscript ``$s$'' stands for ``stable'', and
$\overline{{\mathbb P}^s_0}$ is a complex projective space) together
with a surjective morphism
\begin{equation}\label{eqphi}
\phi\, :\, {\mathbb P}^s_0\, \longrightarrow\, \cN_X
\end{equation}
satisfying the following condition: If ${\mathcal F}\,
\longrightarrow\, X\times T$ is a family, parametrized by $T$,
of stable vector bundles on $X$ of rank $r$ and trivial
determinant, then the corresponding classifying morphism
\[
T \, \longrightarrow\, \cN_X
\]
lifts Zariski locally (in $T$) to maps to ${\mathbb P}^s_0$.
(See \cite[p.\ 81, Sec.\ 7.2]{DN}; the
variety ${\mathbb P}^s_0$ is defined
in \cite[p.\ 86, Remarques 1]{DN}.) In particular, for the
vector bundle $\mathcal E$ in Eqn.\ \eqref{evb} we have
a commutative diagram
\begin{equation}\label{eq.matrix}
\begin{matrix}
{\mathcal U}&\stackrel{\iota}{\hookrightarrow} &
\widehat{\mathbb P}_{x_0}\\
\Big\downarrow && \Big\downarrow f\\
{\mathbb P}^s_0 &\stackrel{\phi}{\longrightarrow} &
\cN_X
\end{matrix}
\end{equation}
where ${\mathcal U}$ is a Zariski open dense subset
of the variety ${\mathcal U}'$
in Eqn.\ \eqref{f1}, and $\phi$ is the map
in Eqn.\ \eqref{eqphi}.

We know that
$\text{Pic}({\mathbb P}^s_0)\, =\, {\mathbb Z}$
(see \cite[p.\ 89, Proposition 7.13]{DN}).
Therefore, the complement
$\overline{{\mathbb P}^s_0}\setminus {\mathbb P}^s_0$
in Eqn.\ \eqref{eqs} is of codimension at least two.
Consequently, from \cite[p.\ 136, Corollaire (6.2)]{Gr},
\begin{equation}\label{eqs2}
\text{Br}'({\mathbb P}^s_0) \, =\, \text{Br}'(
\overline{{\mathbb P}^s_0})\, =\, 0\, .
\end{equation}

Since ${\mathcal U}$ is a nonempty Zariski open subset
of $\widehat{\mathbb P}_{x_0}$, the pull back homomorphism
\[
\iota^*\, :\, \text{Br}'(\widehat{\mathbb P}_{x_0})
\, \longrightarrow\, \text{Br}'({\mathcal U})
\]
is injective, where $\iota$
in the inclusion map in Eqn.\ \eqref{eq.matrix};
see \cite[p.\ 136, Corollaire (6.2)]{Gr}.

Consider the diagram of Brauer groups associated to the diagram in
Eqn.\ \eqref{eq.matrix}. In this diagram, the homomorphism $\iota^*$
is injective, and Eqn.\ \eqref{eqs2} holds. Therefore, we conclude
that Eqn.\ \eqref{f.st2} holds. This completes the proof of the
proposition.
\end{proof}

\begin{remark}
{\rm We note that
Eqn.\ (\ref{eq.Br.}) also follows from \cite[Theorem\ 6.5]{Mu}.}
\end{remark}

By Lemma \ref{lem:1} and the fact that $\cU\subset \cM_X^o\subset 
\cM_X$, we know that
 $$
 \Pic(\cM_X)= \Pic(\cM_X^o)= \Pic(\cU)=\Pic(\cN_X)=\ZZ\, .
 $$
By Proposition \ref{prop:H2}, we know that 
$H^2(\cM_X^o,\ZZ)=\Pic(\cM_X^o)$.
There is a well-defined element
 \begin{equation}
   \label{eqn:extrabla}
   \gamma\in H^2(\cM_X,\ZZ)
 \end{equation}
which is the image of $1$ under
$$\ZZ\,=\,\Pic(\cM_X)\,\longrightarrow\,
H^2(\cM_X,\, \ZZ)\, .
$$
Note that this provides a splitting of the natural map
 $$
 H^2(\cM_X,\ZZ) \surj H^2(\cM_X^o,\ZZ)\, .
 $$

Set $m\, :=\, (r^2-1)(g-1)-3$. Let
 \begin{equation}\label{eqn:hom.prime}
 {F}\, :\, \bigwedge\nolimits^2 H^3(\cM_X,\, \QQ)
 \, \too\, H^{2(r^2-1)(g-1)}(\cM_X,\, \QQ)
 \end{equation}
be the homomorphism defined by $(\alpha\bigwedge\beta)
\,\longmapsto \, \alpha\bigcup\beta\bigcup \g^{\ox m}$.

Consider the exact sequence
 \begin{equation}\label{eqn:extra}
 H^3(\cM_X,\cM_X^o) \,\too\, H^3(\cM_X) \,\stackrel{q}{\too}\, 
H^3(\cM_X^o)\, .
 \end{equation}
We want to show that ${F}$ in Eqn.\ (\ref{eqn:hom.prime}) induces a 
bilinear map on $H^3(\cM_X^o)$. For this we
need the following two lemmas.

\begin{lemma} \label{lem:extra1}
The homomorphism $q\,:\,H^3(\cM_X)\,\too\, H^3(\cM_X^o)$ in Eqn.\ 
(\ref{eqn:extra}) is surjective.
\end{lemma}

\begin{proof}
This map can be defined in families over the moduli space of smooth
projective curves of genus $g$. If it is non-zero for a particular 
curve, then it is non-zero for a generic curve. For a generic curve $X$, 
$H^3(\cM_X^o)\cong H^1(X)$ is
an irreducible Hodge structure. Therefore the map $q$ is surjective for 
such $X$.
But $q$ is a map of local systems over the moduli space of curves. 
Therefore if it is
surjective for a particular point $X$, then it is surjective for all 
curves.

It remains to show that $q$ is non-zero for a particular $X$.
We use the diffeomorphism provided by Eqn.\ (\ref{eqn:extra2}).
Under this diffeomorphism, $\cM_X^o$ corresponds to the smooth
locus $\cR_{g,r}^o$ of $\cR_{g,r}$. So we want to prove
that the homomorphism
 $$
 H^3(\cR_{g,r})\,\too\, H^3(\cR_{g,r}^o)
 $$
is non-zero. By the description in Eqn.\ (\ref{eqn:extra2}), 
$\cR_{g,r}$ is the GIT quotient
$V // {\rm SL}(r, {\mathbb C})$, where $V\,=\,\Hom(\pi_1(X),{\rm SL}(r, 
{\mathbb C}))$.
Therefore $H^3(\cR_{g,r})$ is the ${\rm SL}(r, {\mathbb C})$-invariant 
part of $H^3(V)$. But ${\rm SL}(r, {\mathbb C})$ acts trivially on the 
cohomology, since
it is a connected group. So the pull-back map yields an 
isomorphism 
 $$
 H^3(\cR_{g,r})\cong H^3(V).
 $$

Let $V^o$ be the preimage of $\cR_{g,r}^o$ under the natural map
$V \,\too\, \cR_{g,r}$. We need to check that the homomorphism
$$H^3(V)\cong H^3(\cR_{g,r}) \,\too\, H^3(\cR_{g,r}^o)$$
is non-zero. For this it is enough to check 
that the map
$H^3(V)\,\too\, H^3(V^o)$ is non-zero.

We already know that $$H^3(\cR_{g,r}^o)\,\cong\, H^3(\cM_X^o)\,
\cong\, H^3(\cN_X) \,\cong\, H^1(X)\, .$$
Actually, there is a
universal bundle $\cE \,\too\, X\x \cN_X$ whose second Chern class 
$c_2(\cE)$ produces the required isomorphism
 $$
 c_2(\cE) \ (\, \cdot \,) : H^1(X)
\,\stackrel{\cong}{\too}\, H^3(\cN_X)\, .
 $$ 
We pull back $\cE$ to a bundle over $X\x \cR_{g,r}^o$, and then to a 
bundle $\cE'\,\too\, X\x V^o$.
This bundle is clearly the restriction of the universal bundle 
$\cE''\,\too\, X\x V$.
The bundle $\cE'$ produces a map in cohomology $H^1(X)\,\too\, 
H^3(V^o)$. 
As $\cE'$ extends to the
bundle $\cE''$, we have that the image of the map
$H^3(V)\,\too\, H^3(V^o)$ 
contains the image of the homomorphism
$$H^3(\cR_{g,r}^o) \cong H^1(X)\,\too\, H^3(V^o)\, .$$ 
So we only need to see that the map
$H^3(\cR_{g,r}^o) \,\too\, H^3(V^o)$ is 
non--zero.
But $V^o \,\too\, \cR_{g,r}^o$ is a ${\rm PGL}(r,{\mathbb C})$-bundle, 
so that
$H^3(\cR_{g,r}^o) \,\too\, H^3(V^o)$ is injective. The result is 
complete.
\end{proof}

\begin{lemma} \label{lem:extra2}
The kernel of $F$ defined in Eqn.\ (\ref{eqn:hom.prime})
coincides with the 
image of the homomorphism
$$
H^3(\cM_X,
\cM_X^o)\,\too\, H^3(\cM_X)\, .
$$
Therefore ${F}$ descends to a pairing 
$\overline{F}$ on $H^3(\cM_X^o)\cong H^1(X)$.
\end{lemma}

\begin{proof}
Let $a\in H^3(\cM_X)$ be an element in the image of
$H^3(\cM_X,\cM_X^o)$. We want to check that $a\cup \gamma^{m-3}
\in H^{2m-3}(\cM_X)$ is the zero element.

For this, let us characterize a suitable representative of $\gamma\in 
H^2(\cM_X)$.
By definition, this is the image of the generator of
$$
\Pic(\cM_X)\,=\,
\Pic(\cM_X^o)\,=\,\Pic(\cN_X)\,=\,\ZZ\, .
$$
Take a hyperplane section $H$ of $\cN_X$, consider its preimage 
$\Phi^*H$ in
$\cU \subset \cM_X$ under the projection $\Phi$ of 
Eqn.\ (\ref{eqn:Ph}). Denote by $\overline{H}\subset \cM_X$ the closure 
of this preimage.
Now take $m-3$ generic such choices $\overline{H}_1, \ldots, 
\overline{H}_{m-3}$.
The intersection
 $$
 \overline{W} =\bigcap_{i=1}^{m-3} \overline{H}_i \subset \cM_X
 $$
satisfies that
$\overline{W}\cap \cU=\Phi^*W$, where $W=\bigcap_{i=1}^{m-3} H_i 
\subset \cN_X$.
By Eqn.\ (\ref{cZ}) below, the codimension of the locus of strictly 
semistable bundles $\cZ''$ in the
moduli of semistable bundles $\overline{\cN}_X$ is at least $5$, so the 
intersection of
$m-3$ generic hyperplanes is a projective subvariety not intersecting 
$\cZ''$. Hence $W$
is a $3$-dimensional projective smooth variety. As $\overline{W}$ 
should be connected, it must be that
$\overline{W}\,\subset\, \cU$.
This element represents $\gamma^{m-3}$.

As $a\cup \gamma^{m-3}$ is represented by a cycle which is the 
restriction of $a$ to
$\overline{W}$, and as we may take $a$ supported in a neighborhood of
$\cM_X-\cM_X^o$, hence far away from $\overline{W}$, we have
that $a\cup \gamma^{m-3}\,=\,0$.
\end{proof}

Lemmas \ref{lem:extra1} and \ref{lem:extra2} imply that the map 
in Eqn.\ (\ref{eqn:hom.prime}) induce a map
 \begin{equation}\label{eqn:hom.}
 \overline{F}\, :\, \bigwedge\nolimits^2 H^3(\cM_X^o,\, \QQ)
 \, \too\, H^{2(r^2-1)(g-1)}(\cM_X,\, \QQ)\, .
 \end{equation}

We want to prove that $\overline{F}$ is a polarization for 
$J^2(\cM_X^o)$.

\begin{proposition}\label{prop:4}
The dimension of the image of the homomorphism $\overline{F}$ in
Eqn.\ \eqref{eqn:hom.} is one.
\end{proposition}

\begin{proof}
This is equivalent to prove that the dimension of the image of the 
homomorphism $F$ in
Eqn.\ \eqref{eqn:hom.prime} is one.

We will use the properties the moduli space of Higgs
bundles over $X$ which is naturally homeomorphic to $\cM_X$.
The moduli space of Higgs bundles
is known as the \textit{Dolbeault moduli space}
(see \cite{Si3}).

Let $\cH_X$ denote the moduli space of all semistable
Higgs bundles $(E\, ,\theta)$ over $X$ of the following form:
\begin{itemize}
\item $E$ is a holomorphic vector bundle over $X$ of rank $r$
with $\bigwedge^r E\, \cong\, \cO_X$, and
\item $\theta\, :\, E\, \longrightarrow\, E\bigotimes K_X$
is a Higgs field such that
$\tr (\theta) \, \in\, H^0(X,\, K_X)$ vanishes identically.
\end{itemize}
It is known that the moduli space $\cH_X$ is naturally
homeomorphic to $\cM_X$ \cite[p.\ 38, Theorem 7.18]{Si3}.

Since $\cH_X$ and $\cM_X$ are homeomorphic, 
we have that $H^2(\cH_X,\, \QQ)\, =\, H^2(\cM_X,\, \QQ)$. Let
$\gamma$ be the element of $H^2(\cH_X,\, \QQ)$ corresponding to 
$\gamma\in H^2(\cM_X,\, \QQ)$ defined in Eqn.\ (\ref{eqn:extrabla}). 
Let
 \begin{equation}\label{eqn:Ga}
 \Gamma\, :\, \bigwedge\nolimits^2 H^3(\cH_X, \, \QQ)
 \, \too\, H^{2(r^2-1)(g-1)}(\cH_X,\, \QQ)
 \end{equation}
be the homomorphism defined by $(\alpha\bigwedge
\beta)
\,\longmapsto \, \alpha
\bigcup\beta\bigcup\g^{\ox m}$, where $m=(r^2-1)(g-1)-3$.

Comparing the above homomorphism $\Gamma$ with $F$ defined in
Eqn.\ \eqref{eqn:hom.} we conclude that the following lemma
implies that $\dim \im(F)\, \leq\, 1$.

\begin{lemma}\label{lem:3}
The dimension of the image of the homomorphism $\Gamma$, defined
in Eqn.\ \eqref{eqn:Ga}, is at most one.
\end{lemma}

\begin{proof}
To prove this lemma, we consider the \textit{Hitchin map}
$$
 H\, :\, \cH_X\, \too\, \bigoplus_{i=2}^{r} H^0(X,\, K^{\ox i}_X)
$$
defined by $(E\, ,\theta)\, \longmapsto\, \sum_{i=2}^r
\tr(\theta^i)$ \cite{Hi1}, \cite{Hi2}, \cite[p.\ 20]{Si3}. This
map $H$ is algebraic
and proper \cite{Hi2}, \cite[p.\ 291, Theorem 6.1]{Ni},
\cite[p.\ 22, Theorem 6.11]{Si3}. The fiber
of $H$ over $(0\, , \cdots \, ,0)$ is known as the
\textit{nilpotent cone}.
The nilpotent cone is a finite union of
complete subvarieties of $\cH_X$.
It is in fact a Lagrangian subvariety of $\cH_X$ \cite[p.\
648, Th{\'e}or{\`e}me (0.3)]{La}. Hence each component of the
nilpotent cone is a complete subvariety of dimension
$(r^2-1)(g-1)$.

The moduli space $\cH_X$ is equipped with the following holomorphic
action of $\CC^*$:
 $$
 \lambda \cdot(E\, ,\theta) \,=\, (E\, ,\lambda\cdot\theta)\, ,
 $$
where $\lambda\, \in\, \CC^*$ and $(E\, ,\theta)\, \in\, \cH_X$.
Also, $\CC^*$ acts on $\bigoplus_{i=2}^{r} H^0(X,\, K^{\ox i}_X)$
as
 $$
 \lambda \cdot \sum_{i=2}^r \omega_i\, =\,\sum_{i=2}^r \lambda^i
 \cdot\omega_i\, ,
 $$
where $\omega_i\, \in\, H^0(X,\, K^{\ox i}_X)$. The Hitchin map
$H$ is equivariant with respect to these two actions of $\CC^*$.

We will construct a retraction of $\cH_X$ to a
neighborhood of the nilpotent cone.

Restrict the $\CC^*$--action on $\cH_X$ to the subgroup
$\RR^+\,\subset\, \CC^*$. Consider the map
 $$
 B\,: \,\bigoplus_{i=2}^{r} H^0(X,\, K^{\ox i}_X)
\,\too \,\RR_{\geq 0}\,,
 $$
defined by
 $$
 B(\sum_{i=2}^r \omega_i ) := \sum_{i=2}^r
 || \, \omega_i \, ||^{1/i}\, .
 $$
Clearly $B$ is continuous, proper, and it vanishes
only at the origin. Moreover,
 $$
B(t\cdot\sum_{i=2}^r \omega_i )\,=\, t \cdot B(\sum_{i=2}^r \omega_i )
 $$
for any $t\in \RR^+$. Hence for all $\epsilon\, >\, 0$,
\[
V_\epsilon\,:=\, B^{-1}([0,\epsilon])
\]
is a compact neighborhood of the origin, and
the preimage of
$V_\epsilon$ by the Hitchin map,
 $$
U_\epsilon\,:=\, H^{-1}(V_\epsilon)
$$
is a compact neighborhood of the nilpotent cone in $\cH_X$.
Any open neighborhood of $0\, \in\,
\bigoplus_{i=2}^{r} H^0(X,\, K^{\ox i}_X)$ contains
$V_\epsilon$ whenever $\epsilon$ is sufficiently small.
Since the Hitchin map $H$ is proper, this
implies that any open neighborhood of
$H^{-1}(0)$ contains $U_\epsilon$ provided $\epsilon$ is
sufficiently small.

We have a retraction of $\cH_X$ onto $U_\epsilon$ defined as
follows
$$
\begin{aligned}
R\,:\, \cH_X \x [0\, ,1] &\too\, \cH_X \\
((E,\theta),  t) &\longmapsto \left\{
\begin{array}{ll} (E, t \cdot \theta) , \quad & t\in [0,1]\, ,\ t \geq
\displaystyle\frac{\epsilon}{B(H(E,\theta))}\,, \\[15pt]
(E, t_0 \cdot \theta) , & t \in [0,1]\, ,\ t
\leq t_0 =
\displaystyle\frac{\epsilon}{B(H(E,\theta))} \leq 1 \,,\\[15pt]
(E, \theta) , & t \in [0,1]\, ,\ B(H(E,\theta)) \leq \epsilon\,. \\
\end{array}
\right.
\end{aligned}
$$
Note that in the first two cases, either $t\neq 0$ or $t_0\neq 0$
ensuring that the action is well defined.
For any $(E\, ,\theta)\,\in\, U_\epsilon$, we have
$R((E\, ,\theta),t)\,=\, (E\, ,\theta)$.
Also, $R((E\, ,\theta)\, ,1)\,=\,(E\, ,\theta)$ for
each $(E\, ,\theta)\, \in\, \cH_X$. For all $(E\, ,\theta)\, \in\,
\cH_X$, it can be shown that
$R((E\, ,\theta)\, ,0) \,\in\, U_\epsilon$. Indeed, it is evident
for all $(E\, ,\theta)$ with
$B(H(E,\theta)) \,\leq \,\epsilon$. If $B(H(E\, ,\theta)) \,
\geq \,\epsilon$, then it also holds because
 $$
 B(H(R((E\,,\theta)\, ,0))) \,= \,B(H((E, t_0 \cdot \theta) ))
\,=\, t_0\cdot B(H(E,\theta))\,=\,\epsilon\, .
 $$

Now the nilpotent cone $N=H^{-1}(0)$ is a closed subvariety of
$\cH_X$. Therefore there exists an analytic
open neighborhood $U$ of $N$ (in the smooth topology) that
retracts to $N$, via a retraction $R'$. Take $\epsilon>0$ small
enough so that $U_\epsilon \subset U$. So the above retraction $R$
followed by the retraction $R'$ gives a retraction of $\cH_X$ onto
the nilpotent cone, as required.

The nilpotent cone is a union of components which are complete
subvarieties of (complex) dimension $(r^2-1)(g-1)$.
Therefore, each component of the nilpotent cone
defines an element of $H_{2(r^2-1)(g-1)}(\cH_X,\, \QQ)$.
Using the above retraction it follows that
these elements together generate $H_{2(r^2-1)(g-1)}(\cH_X,\, \QQ)$.

Now the proof of the lemma is completed by the argument
involving monodromy in the proof of Lemma 4.3 in \cite{BM}.
However, one point should be clarified. Unlike in the coprime
case, the monodromy of cohomology ring of moduli spaces of
trivial determinant vector bundles does not factor through
the symplectic group \cite{CLM}. However, we have natural
isomorphisms
\[
H^3(\cM_X,\,{\mathbb Q})\, =\,
H^3(\cN_X,\, {\mathbb Q})\, =\, H^1(X,\, {\mathbb Q})
\]
that extend for a family of curves; see
Proposition \ref{prop:3} and Proposition \ref{prop2}.
Therefore, the monodromy of $H^3(\cM_X,\,{\mathbb Q})$
for a family of curves factors through the symplectic group.
\end{proof}

Continuing with the proof of Proposition \ref{prop:4}, from Lemma
\ref{lem:3} it follows that
 $$
 \dim \im (F)\, \leq\, 1\, .
 $$
We will complete the proof of the proposition by showing that $\im
(F)\, \not=\, 0$.

Let $h$ be a sufficiently positive
class of the smooth quasiprojective variety $\cN_X$
(defined in Section \ref{sec3}). By
Eqn.\ \eqref{eqn:H2-Pic}, the image of $h$ in $H^2(\cN_X,\, \QQ)$,
 \begin{equation}\label{eqn:g1}
 \g\, \in\, H^2(\cN_X,\, \QQ)\,,
 \end{equation}
is a generator of $H^2(\cN_X,\, \QQ)=\QQ$.

Denote by $\overline{\cN}_X$ the moduli
space of (S--equivalence classes of)
semistable vector bundles $E$ on $X$ of rank $r$
and $\bigwedge^r E\,=\,\cO_X$. This moduli space
$\overline{\cN}_X$ is an irreducible projective normal singular
variety with $\cN_X$ as a Zariski open subset.
We will show that the complement
\begin{equation}\label{cZ}
\cZ''\, :=\, \overline{\cN}_X \setminus \cN_X
\end{equation}
is of codimension at least five.

To prove this, for each $\ell\, \in\ [1\, ,r-1]$, let
$\overline{\mathcal N}_X^{\, \ell}$ be the moduli space of semistable
vector bundles $E$ over $X$ of rank $\ell$
and $\bigwedge^\ell E\, =\, {\mathcal O}_X$. Consider the morphism
\[
f_\ell\, :\, \overline{\mathcal N}_X^{\, \ell}\times \overline{\mathcal
N}_X^{\, r-\ell} \times \text{Pic}^0(X)\, \longrightarrow\, \cZ''
\]
defined by $(E_1\, , E_2\, , L)\, \longmapsto\, (E_1\bigotimes
L^{\otimes (r-\ell)})\bigoplus (E_2\otimes (L^*)^{\otimes \ell})$.
We have
\[
\dim \text{image}(f_\ell) \, =\, (\ell^2-1)(g-1)
+((r-\ell)^2-1)(g-1)+g \, \leq\, \dim \cN_X -5\, ,
\]
and also $\cZ''\,=\, \bigcup_{\ell =1}^{r-1} \text{image}(f_\ell)$.
Hence, the codimension of $\cZ''\, \subset\,\overline{\cN}_X$ is at
least five.

Since the codimension of $\cZ''\, \subset\,\overline{\cN}_X$ is at
least five, from \cite[p.\ 55, Th\'eor\`eme B]{DN} and \cite[p.\ 76,
Lemme 5.2]{DN} we know that
$$
\Pic(\overline{\cN}_X)\,=\,\Pic(\cN_X)\,=\,\ZZ\, .
$$
Let $\overline{h}\,\in\,
\Pic(\overline{\cN}_X)$ be the ample class corresponding to $h$. We
assume that $h$ is sufficiently positive so that $\overline{h}$ is
very ample.

Let 
\begin{equation}\label{bagamma}
 \overline{\gamma}\, \in\, H^2(\overline\cN_X,\, \QQ)
\end{equation}
be the cohomology class defined by $\overline{h}\,\in\,
\Pic(\overline{\cN}_X)$. Hence $\overline{\gamma}$ maps to $\gamma$ in 
Eqn.\ \eqref{eqn:g1} by the map 
 $$
 H^2(\overline{\cN}_X,\, \QQ)\, \too \, H^2({\cN}_X,\, \QQ) \, .
 $$
 
Take
\begin{equation}\label{m}
m\,:=\, (r^2-1)(g-1)-3
\end{equation}
(recall that $\dim_\CC\cN_X\, =\,
(r^2-1)(g-1)$). Then a general $(m-1)$--fold intersection
of hyperplanes from the complete linear system
$\vert\overline{h}\vert$ on $\overline{\cN}_X$
does not intersect $\cZ''$, and
furthermore, the intersection is
a smooth projective $4$--fold on $\cN_X$
because the codimension of
$\cZ''$ in Eqn.\ \eqref{cZ} is at least five. Let
\[
S\, \subset\, \cN_X
\]
be a smooth complete $4$--fold obtained by taking $(m-1)$--fold
intersection of hyperplanes from $\vert\overline{h}\vert$.

By the Lefschetz hyperplane section theorem,
\cite[p.\ 215, Lefschetz Theorem (a)]{BSr}, we have
$$H^2(S,\,\QQ) \,=\,H^2(\overline{\cN}_X,\,\QQ)$$
and
\begin{equation}\label{eq-1}
H^3(S,\,\QQ)\,=\, H^3(\overline{\cN}_X,\, \QQ)\, .
\end{equation}

Therefore, using the Hard Lefschetz theorem
for the smooth projective variety $S$, there exists
\begin{equation}\label{eq-4}
\alpha\, , \beta\,\in\, H^3(S,\,\QQ)
\end{equation}
such that the element
$\alpha\bigcup\beta\bigcup \iota^*\overline{\gamma}
\,\in\, H^8(S,\, \QQ)\, =\, {\mathbb Q}$ is
nonzero, where $\iota\, :\, S\, \longrightarrow\,
\overline{\cN}_X$ is the inclusion map, and
$\overline{\gamma}$ is the cohomology class in
Eqn.\ \eqref{bagamma}.
Let $\overline{\alpha}$ (respectively,
$\overline{\beta}$) be the element in
$H^3(\overline{\cN}_X,\, \QQ)$ corresponding to
the element $\alpha$ (respectively,
$\beta$) in Eqn.\ \eqref{eq-4} by the isomorphism
in Eqn.\ \eqref{eq-1}.
We have
\begin{equation}\label{eq-5}
(\overline{\alpha}\cup \overline{\beta}\cup
\overline{\gamma}) \cap [S]\,\neq\, 0
\end{equation}
because $\alpha\bigcup\beta\bigcup
\iota^*\overline{\gamma}$ coincides with
the left--hand side in Eqn.\ \eqref{eq-5}.

Since $S$ is a complete intersection of hyperplanes
on $\overline{\cN}_X$ in the complete linear system
$\vert\overline{h}\vert$ on $\overline{\cN}_X$, and
the first Chern class of $\overline{h}$ is
$\overline{\gamma}$, it follows immediately that the
left--hand side in Eqn.\ \eqref{eq-5}
is a positive multiple of
\[
\overline{\alpha}\cup \overline{\beta}\cup
\overline{\gamma}^{m}\, ,
\]
where $m$ is defined in Eqn.\ \eqref{m}. Therefore, from
Eqn.\ \eqref{eq-5} we have
\begin{equation}\label{eq-6}
\overline{\alpha}\cup \overline{\beta}\cup
\overline{\gamma}^{m}\, \not=\, 0\, .
\end{equation}

Let
$$
\delta\, :\, \overline{\cN}_X\, \longrightarrow\,
\cM_X
$$
be the $C^\infty$
embedding defined by associating to a polystable
vector bundle the unique unitary flat connection on
it (see \cite{NS}). This embedding $\delta$
corresponds to the embedding $\overline{\cN}_X\,
\longrightarrow\, {\mathcal H}_X$ defined by
$E\, \longmapsto\, (E\, , 0)$.

The homomorphism
\begin{equation} \label{eqn:seacaba}
 \delta^*\, :\, H^2(\cM_X,\, {\mathbb Q})\, \longrightarrow\,
 H^2(\overline{\cN}_X,\,
 {\mathbb Q})
\end{equation}
sends $\gamma$ to $\overline{\gamma}$, where 
where $\overline{\gamma}$ is the cohomology
class in Eqn.\ \eqref{bagamma}, and $\gamma$ is defined in Eqn.\ 
\eqref{eqn:extrabla}.

By Lemma \ref{lem:extra1}, the morphism $H^3(\cM_X) \too
H^3(\cM_X^o)\cong H^3(\cN_X)$ is surjective. Let 
 \begin{equation} \label{eqn:queseacaba}
 \alpha,\beta\in H^3(\cN_X)
 \end{equation}
be the images of the classes $\overline{\alpha},\overline{\beta}$ under
$H^3(\overline{\cN}_X)\too H^3(\cN_X)$. 
Take $\widehat{\alpha}, \widehat{\beta} \, \in \, H^3(\cM_X,\QQ)$ some 
preimages of
the cohomology classes in Eqn.\ \eqref{eqn:queseacaba}. Replace
$\overline{\alpha}$ and 
$\overline{\beta}$ by the images of $\widehat{\alpha}$ and 
$\widehat{\beta}$ under 
$\delta^*$ in Eqn.\ (\ref{eqn:seacaba}). Note that
Eqn.\ \eqref{eq-6} continues to hold since it only depends on the 
restriction of 
$\overline{\alpha},\overline{\beta}$ to $S\subset\cN_X$. From
Eqn.\ \eqref{eq-6} it follows immediately that
\[
(\widehat{\alpha}\cup \widehat{\beta}\cup \widehat{\gamma}^m)
\cap [\overline{\cN}_X]
\,=\, \overline{\alpha}\cup \overline{\beta}\cup
\overline{\gamma}^{m}\, \not=\, 0\, .
\]
Consequently, the homomorphism $F$ in Eqn.\ \eqref{eqn:hom.}
is nonzero. We have already shown that
$\dim \im(F)\, \leq\, 1$. Therefore, the proof of the
proposition is complete.
\end{proof}

\section{Torelli type theorems}

In this section we will prove the Torelli theorem for both
$\cM_X$ and $\widehat\cM_X$.

\begin{theorem}\label{thm:1}
Let $X$ and $Y$ be two compact connected Riemann surfaces of genus
$g$, with $g\, \geq\, 3$. Fix $r\, \geq\, 2$. If $g\,=\,3$, then we
assume that $r\,\geq\,3$. Let $\cM_X$ and $\cM_Y$ be the
corresponding moduli spaces of holomorphic ${\rm SL}(r, {\mathbb
C})$--connections defined as in Section \ref{sec:1}. The two
varieties $\cM_X$ and $\cM_Y$ are isomorphic if and only if the two
Riemann surfaces $X$ and $Y$ are isomorphic.
\end{theorem}

\begin{proof}
First note that the smooth locus $\cM_X^o\subset \cM_X$ is determined by
the algebraic structure of $\cM_X$. 
In view of Proposition \ref{prop:3}, to prove the theorem using
the Torelli theorem for curves, \cite{We2}, we need to recover
the principal
polarization on $J^2(\cM_X^o)\, \cong\, {\rm Pic}^0(X)$.

The map $F$ in Eqn.\ \eqref{eqn:hom.prime} is well-defined since the 
element
$\gamma$ is defined by the algebraic structure of $\cM_X$. From
Proposition \ref{prop:4} we know that $\dim \im(F)\, =\, 1$.
Consequently, fixing a generator of $\im(F)$, the homomorphism $F$ 
gives a
nonzero cohomology class
\[
 \theta\, \in\, \bigwedge\nolimits^2 H^3(\cM_X^o,\, \QQ)^* \, =\,
 H^2(J^2(\cM_X^o),\, \QQ)\, .
\]
The one--dimensional subspace
\begin{equation}\label{theta}
\widehat{\theta}\, \subset\, H^2(J^2(\cM_X^o),\, \QQ)
\end{equation}
generated by the above nonzero element $\theta$
is independent of the choices of the generator of $\im(F)$.

Consider the universal curve over the moduli space
${\mathcal M}^0_g$ parametrizing
all smooth complex curves of genus $g$ which do not admit any
nontrivial automorphisms. Let
\[
\rho\, :\, {\rm Pic}^0_{{\mathcal M}^0_g}\, \longrightarrow\,
{\mathcal M}^0_g
\]
be the relative Jacobian for this family. It is known that
the local system $R^2\rho_*{\mathbb Q}$ has exactly one
sub--local system of rank one, and this rank one
sub--local system is generated by the canonical
principal relative
polarization. Indeed, this is an immediate
consequence of the combination of the
fact that the action of $\text{Sp}(2g,
{\mathbb C})$ on $\bigwedge^2 {\mathbb C}^{2g}$ decomposes
it into a direct sum of a one--dimensional
$\text{Sp}(2g, {\mathbb C})$--module and an
irreducible $\text{Sp}(2g, {\mathbb C})$--module
of dimension $g(2g-1)$, and the fact that the monodromy
of $R^1\rho_*{\mathbb C}$ is the full symplectic group
(see \cite[p.\ 710, Theorem 4.2]{BN} for more details).

On the other hand, using the second part
of Proposition \ref{prop:3},
the line $\widehat{\theta}$ in Eqn.\
\eqref{theta} gives a sub--local system of rank one. Therefore,
$\widehat{\theta}$ must be a nonzero rational multiple of the
canonical principal polarization.

Any principal polarization can uniquely be recovered from
any nonzero rational multiple of it; see the proof
of Theorem 4.4 in \cite{BM}.
This completes the proof of the theorem.
\end{proof}

As in Section \ref{sec:1}, let $\widehat\cM_X$ be the moduli
space of holomorphic connections $(E\, ,D)$ with
$\rk(E)\, =\, r$.

\begin{theorem}\label{thm000}
The isomorphism class of the variety $\widehat\cM_X$ determines
the Riemann surface $X$ uniquely up to an isomorphism.
\end{theorem}

\begin{proof}
In view of Theorem \ref{thm:1}, the proof of the theorem is
identical to that of Theorem 5.2 in \cite{BM}. The only point
to note is that there is no nonconstant algebraic map from a
Zariski open subset of ${\mathbb C}{\mathbb P}^1$ to an
abelian variety. The variety $\cM_X$ is unirational.
Indeed, using the diagram in Eqn.\ \eqref{eqn:diag.} it
suffices to show that $\cN_X$ is unirational (recall that
$\Phi$ in Eqn.\ \eqref{eqn:diag.} is an affine bundle).
It is known that $\cN_X$ is unirational; the unirationality
of $\cN_X$
follows from \cite[p.\ 134, Lemma 5.2]{Ne} and \cite[p.\ 136,
Remark]{Ne}.
\end{proof}

\section{Curves defined over real numbers}\label{sec.real}

Fix $t = \sqrt{-1}\,a\, \in\, {\mathbb C}$, where $a$
is a positive real number. Let $\Lambda\, \subset\,
{\mathbb C}$ be the $\mathbb Z$--module generated by
$1$ and $t$. Let
\[
C \,:=\, {\mathbb C}/\Lambda
\]
be the complex elliptic curve.
Consider the anti--holomorphic involution
\[
\sigma\, :\, C\, \longrightarrow\, C
\]
induced by the map ${\mathbb C}\,\longrightarrow\,
{\mathbb C}$ defined by
\[
z \, \longrightarrow\, \overline{z} + \frac{1}{2}\, .
\]
This involution
$\sigma$ of $C$ clearly does not have any fixed points.
Therefore, $(C \,, \sigma)$ is a smooth projective real
curve of genus one without any real points. This real curve
will be denoted by $C_t$.

Consider the real abelian variety
${\rm Pic}^0(C_t)$. We note that there is a canonical
isomorphism
\[
{\rm Pic}^0(C_t)\, \longrightarrow\,
{\rm Pic}^0({\rm Pic}^0(C_t))
\]
that sends any $L\, \in\, {\rm Pic}^0(C_t)$ to the
degree zero divisor
on ${\rm Pic}^0(C_t)$ defined by $L- {\mathcal O}_{C_t}$,
where ${\mathcal O}_{C_t}$ is the trivial line bundle
over $C_t$. Equivalently, the line bundle on ${\rm Pic}^0(C_t)
\times {\rm Pic}^0(C_t)$ defined by the divisor
$\Delta - {\rm Pic}^0(C_t)\times \{{\mathcal O}_{C_t}\}$
gives the identification of ${\rm Pic}^0(C_t)$
with ${\rm Pic}^0({\rm Pic}^0(C_t))$, where $\Delta$ is
the diagonal divisor.

On the other hand, $C_t$ is not isomorphic to
${\rm Pic}^0(C_t)$. Indeed, ${\rm Pic}^0(C_t)$ has a real
point ${\mathcal O}_{C_t}$, while $C_t$ does not have any
real points. Also, note that an abelian variety of dimension
one has exactly one principal polarization.

Therefore, ${\rm Pic}^0(C_t)$ and
${\rm Pic}^0({\rm Pic}^0(C_t))$ are isomorphic
as principally polarized abelian varieties, while
$C_t$ and ${\rm Pic}^0(C_t)$ are not isomorphic.
Consequently, the Torelli theorem fails for real curves
of genus one.

However the following is valid.

\begin{lemma}\label{rem.real}
Let $Y$ and $Z$ be two geometrically irreducible smooth
real projective curves of genus $g_0$, with $g_0\,
\geq\, 2$, such that
\begin{itemize}
\item the real abelian variety ${\rm Pic}^0(Y)$ is
isomorphic to ${\rm Pic}^0(Z)$, and
\item there is an isomorphism ${\rm Pic}^0(Y)
\, \longrightarrow\, {\rm Pic}^0(Z)$ that
takes the canonical principal polarization on
${\rm Pic}^0(Y)$ to that on ${\rm Pic}^0(Z)$.
\end{itemize}
Then the two real algebraic
curves $Y$ and $Z$ are isomorphic.
\end{lemma}

\begin{proof}
Let $Y_{\mathbb C}\, :=\, Y\times_{\mathbb R}
{\mathbb C}$ be the complexification of $Y$.
We know that ${\rm Pic}^0(Y_{\mathbb C})$
is the complexification of the real abelian
variety ${\rm Pic}^0(Y)$.
Furthermore, the canonical
principal polarization on
${\rm Pic}^0(Y_{\mathbb C})$ is given by the
canonical principal polarization on
${\rm Pic}^0(Y)$.

Since the two principally polarized abelian
varieties ${\rm Pic}^0(Y)$ and ${\rm Pic}^0(Z)$
are isomorphic, we know that
${\rm Pic}^0(Y_{\mathbb C})$ is isomorphic, as
a principally polarized variety, to the
Jacobian of $Z_{\mathbb C}\, :=\, Z
\times_{\mathbb R} {\mathbb C}$. Therefore,
the Torelli theorem says that $Y_{\mathbb C}$
is isomorphic to $Z_{\mathbb C}$.

Since $Y_{\mathbb C}$
is isomorphic to $Z_{\mathbb C}$, we may,
and we will, consider $Z$ as a real
structure on the complex curve $Y_{\mathbb C}$. In other words,
the two real curves $Y$ and $Z$ are given by
two anti--holomorphic involutions of
$Y_{\mathbb C}$. Let $\sigma_Y$ (respectively,
$\sigma_Z$) be the anti--holomorphic involutions of
$Y_{\mathbb C}$ defining the real
curve $Y$ (respectively, $Z$).

Let $\tau_Y$ (respectively,
$\tau_Z$) be the anti--holomorphic involutions of
${\rm Pic}^0(Y_{\mathbb C})$ induced by
$\sigma_Y$ (respectively, $\sigma_Z$). Therefore,
$({\rm Pic}^0(Y_{\mathbb C})\, ,\tau_Y)$
and $({\rm Pic}^0(Y_{\mathbb C})\, ,\tau_Z)$
define the real abelian varieties
${\rm Pic}^0(Y)$ and ${\rm Pic}^0(Z)$ respectively.

The anti--holomorphic involution of the Picard
group of the complexification of a
geometrically irreducible smooth projective
real curve preserves the canonical polarization on
it. In particular, both $\tau_Y$ and $\tau_Z$
preserve the element in
$H^2({\rm Pic}^0(Y_{\mathbb C})\, ,{\mathbb Q})$
given by the canonical polarization.

Let
\begin{equation}\label{eta}
\eta\, :\, {\rm Pic}^0(Y_{\mathbb C})\,
\longrightarrow\, {\rm Pic}^0(Y_{\mathbb C})
\end{equation}
be the holomorphic automorphism that satisfies
the identity
\begin{equation}\label{eta2}
\eta\, =\, \tau^{-1}_Z\circ \tau_Y\, ,
\end{equation}
where $\tau_Y$ and $\tau_Z$ are defined above.

Let $C$ be a compact connected complex curve of
genus at least two. The Torelli theorem says that
the group of all automorphisms of the complex Lie group
${\rm Pic}^0(C)$ that preserve the canonical
polarization on it is generated by $\text{Aut}(C)$ together
with the inversion
\[
\iota\, :\, {\rm Pic}^0(C)\, \longrightarrow\,
{\rm Pic}^0(C)
\]
defined by $L\, \longmapsto\, L^*$;
see \cite[p.\ 35, Hauptsatz]{We2}. We note that
the inversion $\iota$ commutes with any automorphism
of the Lie group ${\rm Pic}^0(C)$. Therefore, we
have a homomorphism from $\text{Aut}({\rm Pic}^0(C))$ to
a quotient $Q$ of $\text{Aut}(C)$,
\begin{equation}\label{eta3}
\begin{matrix}
&&\text{Aut}({\rm Pic}^0(C))\\
&&\,{~}\Big\downarrow\phi_C\\
\text{Aut}(C) &\longrightarrow& {Q~\,} &\longrightarrow& 0
\end{matrix}
\end{equation}
The quotient $Q$ is identified with the subgroup of
${\rm Pic}^0(C)$ generated by $\text{Aut}(C)$.

We also know that the natural homomorphism
\begin{equation}\label{def.p}
p\, :\, \text{Aut}(C)\, \longrightarrow\, \text{Aut}(
{\rm Pic}^0(C))
\end{equation}
is injective \cite[p.\ 287, Theorem]{FK}. (The Torelli
theorem fails for real algebraic curves of genus one
because the homomorphism $p$ in that case is not
injective.) Therefore,
the quotient group $Q$ in Eqn.\ \eqref{eta3} actually
coincides with $\text{Aut}(C)$. Furthermore, the
homomorphism $\phi_C$ in Eqn.\ \eqref{eta3} satisfies
the identity
\begin{equation}\label{eta0}
\phi_C\circ p\, =\, \text{Id}_{\text{Aut}(C)}\, ,
\end{equation}
where $p$ is the homomorphism in Eqn.\ \eqref{def.p}.

Define $\phi_Y$ as in Eqn.\ \eqref{eta3} by substituting
the curve $Y$ for $C$. Consider $\eta$ constructed in
Eqn.\ \eqref{eta}. Since it satisfies
Eqn.\ \eqref{eta2}, using Eqn.\ \eqref{eta0}
we conclude that the automorphism
\[
\phi_Y(\eta)\, :\, Y\, \longrightarrow\, Y
\]
satisfies the identity $\phi_Y(\eta)
\, =\, \sigma^{-1}_Z\circ \sigma_Y$.

In other words, $\phi_Y(\eta)$ is an isomorphism
between the two real curves $(Y_{\mathbb C}\, ,
\sigma_Y)$ and $(Y_{\mathbb C}\, , \sigma_Z)$.
This completes the proof of the lemma.
\end{proof}

In view of Lemma \ref{rem.real}, we conclude that
the Torelli theorems in \cite{MN} and \cite{NR}
for the moduli space of stable vector bundles remain valid
for curves defined over $\mathbb R$. Consequently,
Theorem 4.4 and Theorem 5.2 of \cite{BM} remain
valid for curves defined over $\mathbb R$.

Similarly, \cite[p.\ 2, Theorem 1.0.2]{AS} remains valid for
curves defined over $\mathbb R$. Therefore, Theorem
\ref{thm:1} and Theorem \ref{thm000} proved here also
remain valid for curves defined over $\mathbb R$.

\bigskip

\noindent \textbf{Acknowledgements.}\, We thank Tam\'as
Hausel, Frances Kirwan, Pierre Deligne, Lisa Jeffrey and
Peter Newstead for useful comments. We thank the referee
for a careful reading of the manuscript. The second author wishes to
acknowledge Universidad Complutense de Madrid and
Institute for Advanced Study at Princeton for their
hospitality and for providing excellent working conditions.


\end{document}